\documentclass[12pt]{amsart}
\usepackage{amsmath, amssymb, euscript, hyperref}

\newtheorem{thm}{Theorem}[section]
\newtheorem{lem}[thm]{Lemma}

\theoremstyle{definition}
\newtheorem{defn}[thm]{Definition}

\newtheorem{rem}[thm]{Remark}
\newtheorem{cor}[thm]{Corollary}

\address{Azer Akhmedov, Department of Mathematics,
North Dakota State University,
Fargo, ND, 58108, USA}
\email{azer.akhmedov@ndsu.edu}

\begin{document}

 \begin{center} {\bf {\LARGE Girth Alternative for Subgroups of $PL_o(I)$}} \end{center}

\vspace{1cm}

 \begin{center} Azer Akhmedov \end{center}

 \bigskip

\begin{center} {\small ABSTRACT: We prove the Girth Alternative for finitely generated subgroups of $PL_o(I)$. We also prove that a finitely generated subgroup of {\em Homeo}$_{+}(I)$ which is sufficiently rich with hyperbolic-like elements has infinite girth.} \end{center}

\vspace{1cm}

 \section{Introduction} 
 
 \medskip
 
 The notion of a girth for a finitely generated group was first introduced in [S] motivated by the study of Heegaard splittings of closed 3-manifolds. 

 \medskip

 \begin{defn} Let $\Gamma $ be a finitely generated group. For any finite generating set $S$ of $\Gamma $, $girth (\Gamma , S)$ will denote the minimal length of relations among the elements of $S$. Then we define $girth (\Gamma ) = \displaystyle \sup _{\langle S \rangle = \Gamma , |S| < \infty }girth (\Gamma , S)$. 
 \end{defn}

\medskip

 By definition above, an infinite cyclic group has infinite girth, but this fact should be viewed as a degeneracy since (as remarked in [Ak1]) any group satisfying a law and non-isomorphic to $\mathbb{Z}$ has a finite girth. 

 \medskip

 In [Ak2], we have proved that if a finitely generated group is word hyperbolic, or one-relator, or linear then it is either virtually solvable or has infinite girth.  More generally, given a class $\mathcal{C}$ of finitely generated groups, we will say that $\mathcal{C}$ satisfies Girth Alternative if any group from the class $\mathcal{C}$ is either virtually solvable or has infinite girth. 

\medskip

 In [Y], S.Yamagata has proved the Girth Alternative for convergence groups and for irreducible subgroups of mapping class groups. Independently, in [N1] and [N2], K.Nakamura proves the alternative for convergence groups but also for all subgroups of mapping class groups as well as for subgroups of $Out(\mathbb{F}_n)$ that contain irreducible elements with irreducible powers.

 \medskip

 In this paper, we will prove that the Girth Alternative holds for subgroups of $PL_o(I)$ - the group of orientation preserving piecewise linear homeomorphisms of the closed interval $I = [0,1]$. It is known that any virtually solvable subgroup of $PL_o(I)$ is indeed solvable (see [Bl], Corollary 1.3.) so the Girth Alternative in this case is equivalent to the following

\medskip
 
 \begin{thm} Any finitely generated subgroup of $PL_o(I)$ is either solvable or has infinite girth.
\end{thm}

\medskip

 It is remarked in [Ak1] that a finitely generated non-cyclic group which satisfies a law has a finite girth. Thus we obtain the following corollary.

 \begin{cor} If a subgroup of $PL_o(I)$ satisfies a law then it is solvable.
 \end{cor} 

 \medskip
 
 \begin{rem} As another corollary of Theorem 1.2, we obtain that $girth (F) = \infty $ where $F$ denotes the R.Thompson's group. This fact has been proved in [Br1] and in [AST]; both proofs use different ideas from each other and from the proof of Theorem 1.2. Theorem 1.2 also implies that $girth(B) = \infty $ where $B$ is the Brin group introduced in [Br2] under the notation $G_1$ (the notation $B$ is used in [Bl] and in [DS], among other places).  
\end{rem}

\medskip

 It is easy to prove the Girth Alternative for {\em Diff}$_{\omega }(I)$ - the group of orientation preserving analytic diffeomorphisms of $I$, however, we do not know if the alternative still holds when the regularity is decreased. The following questions are interesting to us:

\medskip

 {\bf Question 1.} Does Girth Alternative hold for subgroups  

\medskip

 a) of {\em Homeo}$_{+}(I)$?   \ \ \ \ \ b) of  {\em Diff}$_{+}(I)$?

\medskip

 {\bf Question 2.} Is there a finitely generated subgroup of $PL_o(I)$ which is not embeddable into  {\em Diff}$_{+}(I)$?

\medskip

 In regard to Question 1, we prove the following partial result. 

\bigskip

\begin{thm} Let $\Gamma $ be any finitely generated subgroup of {\em Homeo}$_{+}(I)$. Assume that for all $N\in \mathbb{N}$, for every sequence $0 < x_1 < x_2 < \ldots < x_N < 1$, and for all $\epsilon > 0$, one can find an element $\gamma \in \Gamma $ such that $Fix(\gamma ) = \{0, c_1, \ldots , c_N, 1\}$, and $|c_i-x_i| < \epsilon $, for all $1\leq i\leq N$.  Then $girth (\Gamma ) = \infty $.
\end{thm}

\medskip

 \begin{rem} As a corollary of Theorem 1.5, we obtain yet another proof of the fact that $girth (F) = \infty $.
\end{rem}

\bigskip

\section{Towers and Exemplary Towers: Review of C.Bleak's Results}

 \medskip
 
 In the proof of Theorem 1.2, as a crucial tool, we use the result of C.Bleak on the existence of arbitrarily high towers in a non-solvable subgroup of $PL_o(I)$ ([Bl]). First, we would like to introduce the following notions essentially borrowed from [Bl], with a slightly different terminology.

 \medskip

 \begin{defn} An ordered $n$-tuple $(f_1, \ldots , f_n)$ of elements of  $PL_o(I)$ is said to form {\em a strict tower} if there exist intervals $(a_i, b_i), 1\leq i\leq n$ such that 

 \medskip

 (i) $0 < a_1 < \ldots < a_n < b_n < \ldots < b_1 < 1$;

\medskip

 \medskip

 (ii) for all $i\in \{1, \ldots , n\}, \ f_i(a_i) = a_i, f(b_i) = b_i$, and $f_i$ has no fixed points in $(a_i, b_i)$.

 \medskip

 (iii) for all $i, j\in \{1, \ldots , n\}$, \ if $i < j$ then $f_i(x) > f_j(x), \forall x\in [a_j,b_j]$

 \medskip

 We will denote the strict tower by $T = [(f_1, \ldots , f_n); (a_1,b_1), \ldots , (a_n, b_n)]$; $n$ will be called {\em the height of the tower $T$}, and the interval $(a_i, b_i)$ will be called {\em the $i$-th base} of the tower. 
 \end{defn}

 \medskip

  \begin{defn}We will say that a strict tower $$T = [(f_1, \ldots , f_n); (a_1,b_1), \ldots , (a_n, b_n)]$$ is {\em suitable } if for any nonzero integer $p$ and for all $1\leq i < j \leq n$, the following condition holds: $ f_i^p([a_j,b_j])\cap \displaystyle \mathop {\cup}_{i+1\leq k\leq n} supp (f_k) = \emptyset  \ (1)$.
  \end{defn}
	
 \bigskip

 \begin{rem} Condition (1) implies that for any nonzero integer $p$ and for all $1\leq i < j \leq n$, $ f_i^p([a_j,b_j])\cap [a_j,b_j] = \emptyset \ (2)$. Notice that, if the $n$-tuple $(f_1, \ldots , f_n)$ of elements of  $PL_o(I)$ forms a tower then for sufficiently big $q\in \mathbb{N}$, the $n$-tuple $(f_1^q, \ldots , f_n^q)$ forms a  tower with the same bases which satisfies condition $(2)$. Also, the existence of a suitable tower of arbitrary height in non-solvable subgroups of $PL_o(I)$ immediately follows from the existence of the {\em exemplary towers} of arbitrary height, [Bl]. 
 \end{rem}
 
 \bigskip
 
 To explain the existence of exemplary towers, we would like to make a digression into some of the results of C.Bleak. For the rest of this section let $\Gamma \leq PL_o(I)$. The following notions and notations are all borrowed directly from [Bl]. 
 
 \medskip
 
 We will call the convex hull of a point in $I$ under the action of $\Gamma $ an {\em orbital} of $\Gamma $ if this convex hull contains more than one point. We note that the orbitals are open intervals. If $g\in \Gamma $, we will refer to an orbital of the group $\langle g \rangle $ as an orbital of $g$. If an open interval $A$ is an orbital of $g$ then the pair $(A,g)$ will be called a {\em signed orbital of} $G$. $g$ will be called the {\em signature} of the signed orbital $(A,g)$.
 
 \medskip
 
 Given a set $Y$ of signed orbitals of $G$, the symbol $S_Y$ will refer to the set of signatures
of the signed orbitals in $Y$. Similarly, the symbol $O_Y$ will refer to the set of orbitals of the
signed orbitals of $Y$. We note that the set of signed orbitals of $PL_o(I)$ is a partially ordered
set under the lexicographic order induced from the partial order on subsets of $I$ (induced by
inclusion) in the first coordinate, and the left total order of the elements of $PL_o(I)$ in the
second coordinate.
  
 \medskip
 
 A {\em tower} $T$ of $G$ is a set of signed orbitals which satisfies the following two criteria.

 \medskip
 
 1. $T$ is a chain in the partial order on the signed orbitals of $G$.
 
 \medskip
 
 2. For any $A\in O_T$, $T$ has exactly one element of the form $(A,g)$.
 
 \medskip
 
  Given a tower $T$ of $G$, if $(A,g), (B,h)\in T$ then one of $A\subseteq B$ and $B\subseteq A$ holds, with equality occurring only if $g = h$ as well. Therefore, one can visualize the tower as a stack of nested levels that are always getting wider as one goes ``up" the stack.
  
 \medskip
 
 The cardinality of the set $O_T$ will be called the {\em height} of the tower $T$.

 \medskip
 
 A major result of [Bl] is the following beautiful geometric characterization of solvable subgroups of $PL_o(I)$.
 
 \begin{thm}If $G\leq PL_o(I)$ is a non-solvable subgroup if and only if $G$ admits a tower of height $n$ for any $n\geq 1$. 
 \end{thm}
 
 \bigskip
 
 If $G$ admits two signed orbitals $(A,g)$ and $(B,h)$ so that $A = (a_1,a_2)$ and $B = (b_1, b_2)$,
with $a_1 < b_1 < a_2 < b_2$ then we will say that $G$ admits a {\em transition chain of length two}.
One can similarly define transition chains of any (finite) length, but we will have no need for that generality here.

 \medskip
 
If $A = (a,b)$ is an orbital of the group $G$, and $G$ has an element $g$ which has an orbital
$B = (c,d)$ so that either $c = a$ or $d = b$, then we say that $g$ has an orbital that shares an
end with $A$.
 
 \medskip
 
 Given an orbital $A$ of $H$ we say that $h$ {\em realizes an end of} $A$ if some orbital of $h$ lies
entirely in $A$ and shares an end with $A$. If $g$ and $h$ are elements of $PL_o(I)$ and there is an interval $B = (a,b)\subset I$ so that both $g$ and $h$ have $B$ as an orbital, then we will say that $g$ and $h$ {\em share the orbital} $B$.
 
 \bigskip
 
 We will say that an orbital $A$ of a group $H\leq PL_o(I)$ is {\em imbalanced} if some element
$h\in H$ realizes one end of $A$, but not the other, and we will say $A$ is {\em balanced} if whenever
an element $h\in H$ realizes one end of $A$, then $h$ also realizes the other end of $A$ (note that
$h$ might do this with two distinct orbitals). A subgroup $H\leq PL_o(I)$ will be called balanced if given any subgroup $G\leq H$, and any orbital $A$ of $G$, every element of $G$ which realizes one end of $A$ also realizes the other end of $A$. In the case where $H$ has a subgroup $G$ with an imbalanced orbital, then we will say that $H$ is imbalanced.

 \medskip
 
 We say a tower $T$ is an {\em exemplary tower} if the
following two additional properties hold:

1. Whenever $(A,g), (B,h)$ implies the orbitals of $g$ are disjoint
from both ends of the orbital $B$.

2. Whenever $(A, g), (B, h)\in T$ then $(A,g), (B,h)$ implies no orbital of $g$ in $B$ shares
an end with $B$.

 \medskip

 C.Bleak proves several technical results which indicate the plethora of exemplary towers in $PL_o(I)$.
The following two lemmas are stated in [Bl] as Lemma 1.4. and Lemma 2.12, respectively. 

 \medskip
 
 \begin{lem} If $H$ is a subgroup of $PL_o(I)$, and $H$ admits a transition chain of length two, then $H$ admits infinite towers. 
 \end{lem}
 
 \medskip
  
 \begin{lem} If $H$ is a subgroup of $PL_o(I)$ which does not admit a transition chain of length
two and $G$ has a tower $T$, then $T$ is exemplary.
\end{lem}
 
 \medskip
 
  Thus, the existence of a transition chain of length two implies the existence of infinite towers, and the absence of the transition chain of length two implies that all towers are nice, i.e. exemplary. It turns out the absence of a transition chain of length two also implies certain nice properties of the group itself. The following result is stated as Remark 4.9. in [Bl]. 
  
  \medskip 
 
 \begin{lem} If $G$ is a subgroup of $PL_o(I)$ that does not admit transition chains of length
two, then $G$ is balanced.
\end{lem}

 \medskip
   
 We also need the following two technical results from [Bl].

 \begin{lem}[Corollary 2.13.(Bl)] If $G$ is a balanced subgroup of $PL_o(I)$ and $G$ admits a tall tower in some orbital $A$, or $G$ admits a deep tower in some orbital $A$, or $G$ admits a bi-infinite tower in
some orbital $A$, then $G$ admits an exemplary tall tower in $A$, or $G$ admits an exemplary deep
tower in $A$, or $G$ admits an exemplary bi-infinite tower in $A$, respectively. 
\end{lem}

 \medskip

 \begin{lem}[Lemma 2.8.(Bl)] Suppose $H\leq PL_o(I)$, and that $G\leq H$ has imbalanced orbital $A = (a,b)$. Then $H$ admits an exemplary bi-infinite tower $E$ whose orbitals are all in $A$.
 \end{lem}
 
 \medskip

 \begin{cor} If $H$ is an imbalanced subgroup of $PL_o(I)$ then $H$ admits an exemplary bi-infinite tower. 
 \end{cor}
 
 \medskip
 
 Combining the above results we can now claim the following lemma.
 
 \medskip
 
 \begin{lem} If $\Gamma \leq PL_o(I)$ is not solvable then it contains an exemplary tower of any height $n\geq 1$. 
 \end{lem}
 
 {\bf Proof.} If $\Gamma $ does not admit a transition chain of length two then the claim follows from Theorem 2.4. and Lemma 2.6. If $\Gamma $ is not balanced then the claim follows from Corollary 2.10. Thus we can assume that $\Gamma $ is balanced and admits a transition chain of length two. Then by Lemma 2.5. $\Gamma $ admits an infinite tower. Then by Lemma 2.8. it admits an exemplary bi-infinite tower. $\square $

 \bigskip
 
 \section {Proof of Theorem 1.2.} 
 
 \medskip
 
 For any natural number $r$, let $G_r = (\ldots ((\mathbb{Z}\wr \mathbb{Z})\wr \mathbb{Z})\wr \ldots \mathbb{Z})\wr \mathbb{Z}$ where the iterated wreath product is taken $r$ times. The group $G_r$ can be defined inductively as $G_0 = 1, G_{i+1} = G_i\wr \mathbb{Z}, 0\leq i\leq r-1$. In the wreath product $G_i\wr \mathbb{Z}$, the standard generator of the acting group $\mathbb{Z}$ will be denoted by $g_{r-i}$. (in [Bl], the group $G_r$ is denoted by $W_r$).

\medskip 

The following lemma will be useful; the idea of its proof is essentially borrowed from  the proof of Lemma 2.3. in [Ak1].

 \bigskip

 \begin{lem} For all $q, k\in \mathbb{N}$, there exist $r\in \mathbb{N}$ and $w_1, \ldots , w_k\in G_r$ such that there is no relation of length less than $q$ among $w_1, \ldots , w_k$.  \ $\square $  
\end{lem}

 \medskip
 
 {\bf Proof.} Since a free group on 2 generators contains a free group on $k$ generators for any $k\geq 3$, it is sufficient to prove the claim for $k=2$. We will do this by induction; more precisely, it suffices to prove the following claim: {\em if $A, B$ are non-trivial groups and $A$ satisfies no law in two variables of length less than $n\geq 4$ then the wreath product $B\wr A = A\ltimes \displaystyle \mathop{\oplus }_{i\in A}B $ satisfies no law in two variables of length less than $n+1$}.
 
 \medskip
 
 Indeed, let $w_1, w_2\in A$ with no relation of length less than $n$. An element of $\displaystyle \mathop{\oplus }_{i\in A}B$ can be written as $(x_g)_{g\in A}$ where all but finitely many ``coordinates" are 1. Let $b$ be a non-identity element of $B$, and $t = (x_g)_{g\in A}$ be the element of $\displaystyle \mathop{\oplus }_{i\in A}B$ where $x_g = b$ for $g=1$ and $x_g = 1$ otherwise.
 
 \medskip
 
 Then there is no relation of length less than $n+1$ among the elements $tw_1$ and $w_2$. Indeed, let $W(x,y)$ be a non-trivial reduced word of length $k < n+1$ in the alphabet $\{x^{\pm 1}, y^{\pm 1}\}$ such that $W(tw_1,w_2) = 1\in B\wr A$. For every $1\leq i\leq k$, let also $W_i$ be the prefix of $W$ of length $i$. Then $W(w_1,w_2) = W_k(w_1,w_2) = 1\in A$, moreover, for at least one $i\in \{1,\ldots , k\}$ there exists $j\in \{1,\ldots, k\}\backslash \{i\}$ such that $W_i(w_1,w_2) = W_j(w_1, w_2)$. But this implies a relation of length less than $n$ among $w_1$ and $w_2$. $\square $
  
\bigskip

 We also would like to observe the following simple lemma.

\medskip

 \begin{lem} Let $T = [(f_1, \ldots , f_n, f_{n+1}); (a_1,b_1), \ldots , (a_n, b_n), (a_{n+1}, b_{n+1})]$ be a suitable tower, and $\phi _1, \ldots , \phi _n\in PL_o(I)$ be such that for all $i\in \{1, \ldots , n\}$, $supp(\phi _i) \subseteq [a_i, b_i]$ and $\phi _i|_{[a_i,b_i]} = f_i|_{[a_i,b_i]}$. Then

\medskip

 a) the maps $\phi _1, \ldots , \phi _n$ generate a subgroup $H_n\leq PL_o(I)$ such that there exists an isomorphism $\Phi : H_n\rightarrow G_n$ given by $\Phi (\phi _i) = g_i, 1\leq i\leq n$; 

\medskip

 b) if $W = W(g_1, \ldots , g_n)$ is any non-identity element of $G_n$ then there exists a word $U = U(g_1, \ldots , g_n)\in G_n$ such that $W_1((U_1(a_{n+1}, b_{n+1})))\cap U_1((a_{n+1}, b_{n+1})) = \emptyset $ where $W_1 = W(\phi _1, \ldots , \phi _n), U_1 = U(\phi _1, \ldots , \phi _n)$; 

\medskip

 c) if $W = W(g_1, \ldots , g_n)\in G_n$ and $x\in (a_{n+1}, b_{n+1})$ then $W(f_1, \ldots , f_n)(x) = W(\phi _1, \ldots , \phi _n)(x)$. \ $\square $ 
\end{lem}

\bigskip

 Now, let $\Gamma $ be a finitely generated non-solvable subgroup of $PL_o(I)$, $s = d(\Gamma )$ (i.e. $s$ is the minimal cardinality for a generating subset of $\Gamma $). Without loss of generality we may assume that $\Gamma $ is irreducible, i.e. it has no global fixed point in the interval $(0,1)$.

\medskip

 For every finite generating set $X$ of $\Gamma $ we will fix the left-invariant Cayley metric on $\Gamma $ with respect to $X$, and let $B_k(1; X)$ denote the ball of radius $k$ centered at identity element $1\in \Gamma $, for all $k\geq 1$. We also let $$S_k(\Gamma , X, c) = \{\gamma \in B_k(1; X) \ | \ \gamma '(0) = c\}, \   S_k(\Gamma , X) =  \{\gamma \in B_k(1; X) \ | \ \gamma '(0) = 1\},$$ \ $$C_k(\Gamma , X) = \{c\in \mathbb{R}_{+} \ | \ S_k(\Gamma , X, c) \neq \emptyset \}.$$     

\medskip

 Fix a positive integer $m$, and let $q = 2m^2$. By Lemma 3.1, there exists $r\in \mathbb{N}$ such that in the group $G_r$ there exist $s$ elements $w_1, w_2, \ldots , w_s$ such that there is no relation of length less than $q$ among $w_1, \ldots , w_s$. Let $g_1, \ldots , g_r$ be the standard generators  of $G_r$ and let $w_i = W_i(g_1, \ldots , g_r)$, $1\leq i\leq s$ where $W_i$ is a reduced word in the free group of rank $r$ formally generated by the letters $g_1, \ldots , g_r$.

\medskip

 Since $\Gamma $ is non-solvable, the commutator subgroup $[\Gamma , \Gamma ]$ is non-solvable. Then by Remark 2.3 and Lemma 2.11, there exists an ordered $(r+1)$-tuple $(f_1, \ldots , f_r, f_{r+1})$ of elements of $[\Gamma , \Gamma ]$ which form a suitable tower of height $r+1$.  Let $0 < d < D < 1$ such that $\displaystyle \mathop{\cup }_{1\leq i\leq r+1}supp(f_i)\subset (d,D)$. 

\medskip

 Then we can find $\epsilon _0 > 0$, a finite generating set $S$ of $\Gamma $ of cardinality $s$, and a suitable tower $(h_1, \ldots , h_{r}, h_{r+1})$ of elements of $\Gamma $ with bases $(a_i, b_i)\subset (0, \epsilon _0), 1\leq i\leq r+1$ such that the following conditions hold:

 \medskip

 (i) for all $\beta \in B_{2m}(1; S), \beta (\Omega )\subset (0, \epsilon _0)$;

\medskip

(ii) for any two distinct $c_1, c_2 \in \mathbb{R}_{+}$, and for all $\beta _1\in S_{2m}(\Gamma , S, c_1), \beta _2\in S_{2m}(\Gamma , S, c_2)$, \ $\beta _1(\Omega )\cap \beta _2(\Omega ) = \emptyset $;

 \medskip

 (iii)  for all $c\in \mathbb{R}_{+}, \beta _1, \beta _2\in S_{2m}(\Gamma , S, c)$, and for all $x\in \Omega $, $\beta _1(x) = \beta _2(x)$. [so, in particular, for all $\beta \in S_{2m}(\Gamma , S)$ and for all $x\in \Omega , \ \beta (x) = x$]; 

 \medskip

 where $\Omega  = \displaystyle \mathop{\cup }_{1\leq i\leq r+1}supp(h_i)$.
  
  \bigskip
  
  Indeed, let $X_0 = \{\alpha _1, \dots , \alpha _s\}$ be any generating set of $\Gamma $ of cardinality $s$. Without loss of generality, we may assume that $$1\leq \alpha _1'(0) \leq \dots \leq \alpha _s'(0) \ \mathrm{and} \ \alpha _s'(0) > 1.$$ 
  
  \medskip
  
  Let $\delta > 0$ such that $\alpha _s$ has no singularity in $(0,\delta )$. By irreducibility, there exists $\phi \in \Gamma $ such that $\phi (D) < \delta $. For $n\geq 1$ and $1\leq i\leq r+1$, we let $f_i^{(n)} = \psi_n^{-1} f_i\psi _n$, where $\psi_n = \phi ^{-1}\alpha _s^n$.  Let also $(d_n, D_n) = \psi _n^{-1}((d,D)), n\geq 1$. Notice that the interval $(d_n, D_n)$ (in particular, its subset  $\displaystyle \mathop{\cup }_{1\leq i\leq r+1}supp(f_i^{(n)})$) converges to zero as $n\to \infty $, moreover, there exists a positive integer $p$ such that for all $k\geq p$, we have $\frac{D_k}{d_k} = \frac{D_p}{d_p}$.
  
  \medskip
  
  We can choose a finite generating set $S = \{\gamma _1, \dots , \gamma _s\}$ such that $\min \{c\in C_{2m}(\Gamma , S) \ | \ c > 1\} > \frac{D_p}{d_p}$. \footnote{First, we let $X_1 = \{\beta _1, \dots , \beta _s\}$ where $\beta _s = \gamma _s$, and $\beta _i = \gamma _i\gamma _{s}^{n_i}, 1\leq i\leq s-1$ for some $n_{s-1}\leq \dots \leq n_2\leq n_1$ such that $1 < \beta _s'(0) \leq \beta _{s-1}'(0) \leq \dots \leq \beta _1'(0)$ and $\beta _{1}'(0) > \frac{D_p}{d_p}$. After this, we modify the generating set $X_1$ further by letting $S = \{\beta _1, \beta _1^{\lambda }\beta _2, \beta _1^{\lambda ^2}\beta _{3}, \dots , \beta _1^{\lambda ^{s-1}}\beta _{s}\}$ where $\lambda  = 4m+1$.}

  \medskip
  
  Then, for sufficiently big $n$, we can take  $h_i = f_i^{(n)}, 1\leq i\leq r+1$ to satisfy the claims (i)-(iii). 
   
	\bigskip

	 Now, let $A_m$ be a minimal subset of $B_m(1;S)$ such that $A_m\cap S_m(\Gamma , S, c) \neq \emptyset $ for every $c\in C_m(\Gamma , S)$; and for all $k\in \{1, \ldots , s\}$, let  
	
	$$\omega _k = \displaystyle \prod _{\gamma \in A_m}(\gamma v_k\gamma ^{-1}), \ \eta _k = \displaystyle \prod _{\gamma \in A_m}(\gamma u_k\gamma ^{-1})$$
	
	\medskip
	
	where $v_k = W_k(h_1, \ldots , h_r), u_k = v_k^m, 1\leq k\leq s$ and for the products (the formulas for $\omega _k$ and $\eta _k$) we choose any linear order on the set $A_m$.

	\medskip
	
	Notice that because of conditions (i)-(iii), for any two $\gamma ', \gamma ''\in B_m(1)$, we have $[\gamma '\omega _k(\gamma ')^{-1}, \gamma ''\omega _k(\gamma '')^{-1}] = 1$ and $[\gamma '\eta _k(\gamma ')^{-1}, \gamma ''\eta _k(\gamma '')^{-1}] = 1$, so the order in the products $\displaystyle \prod _{\gamma \in A_m}(\gamma v_k\gamma ^{-1})$ and $\displaystyle \prod _{\gamma \in A_m}(\gamma u_k\gamma ^{-1})$ does not matter.
	
	\medskip
	
	Now, let $S^{(m)} = \{\eta _1\gamma _1\eta _2, \eta _2\gamma _2 \eta _3, \ldots , \eta _s\gamma _s\eta _1, \omega _1, \ldots , \omega _s\}$. Then $S^{(m)}$ generates $\Gamma $, and there is no relation of length less than $m$ among the elements of $S^{(m)}$. 
	
	\medskip
	
	Indeed, let $R = R(\eta _1\gamma _1\eta _2, \eta _2\gamma _2 \eta _3, \ldots , \eta _s\gamma _s\eta _1, \omega _1, \ldots , \omega _s)$ denotes such a relation. Then we can write $$R(\eta _1\gamma _1\eta _2, \eta _2\gamma _2 \eta _3, \ldots , \eta _s\gamma _s\eta _1, \omega _1, \ldots , \omega _s) = R_0(\theta _1R_1\theta _1^{-1})(\theta _2R_2\theta _2^{-1})\ldots  (\theta _nR_n\theta _n^{-1})R_{n+1}$$

 where $n\leq m$, $\theta _i\in B_m(1)$ for all $1\leq i\leq n$, and $R_j = R_j(\eta _1, \ldots , \eta _s, \omega _1, \ldots , \omega _s)$ is a reduced word of length at most $m$ for all $0\leq j\leq n+1$; moreover, $R_1, \ldots , R_n$ are nontrivial.  

	\medskip
	
	Then we obtain a nontrivial relation $V(v_1, \ldots , v_s)$ among $v_1, \ldots , v_s$ of length at most $2m^2$. We can write $V(v_1, \ldots , v_s) = W(h_1, \ldots, h_r) = 1$ where $W$ is a nontrivial reduced word. Notice that $V(v_1, \ldots , v_s) = W(h_1, \ldots, h_r)$ represents a map in $PL_o(I)$, while $V(w_1,\ldots , w_s)$ represents a word in $G_r$ which by our choices, does not represent an identity element in $G_r$. Then, by Lemma 3.2, for some word $U = U(h_1, \ldots , h_r)$ and for all $x\in U((a_{r+1}, b_{r+1}))$ we have $W(x)\notin U((a_{r+1}, b_{r+1}))$ thus $W(x)\neq x$. Since $m$ is arbitrary, we conclude that $girth(\Gamma ) = \infty $. $\square $
	
 \medskip
 
 The main idea of the proof of Theorem 1.2 is described below. For any given $q, k$, in a suitable tower of sufficiently big height $r$, formed by PL-maps $\phi _1, \ldots , \phi _r, \phi _{r+1}$, one can find words $w_1, \ldots , w_k$ in the alphabet $\{\phi _1^{\pm 1}, \ldots , \phi _r^{\pm 1}\}$ such that the corresponding elements (let us denote them by $\overline{w_1}, \ldots , \overline{w_k}$) do not have a relation of length less than $q$. This is because, upon the action on the innermost base of the tower (an orbital of $\phi _{r+1}$), a suitable tower behaves as if the maps generate a genuine copy of $G_r$, for elements in the ball of certain radius. The problem is how to find a finite generating set without a short relation among the generators, not just among {\em some} elements. For this, we pick up a tower with sufficient height such that the PL homeomorphisms forming this tower are supported in a very small interval. This interval can be made arbitrarily small, therefore, using irreducibility of $\Gamma $, one can push this support (interval) close enough to the end 0 of the interval $[0,1]$ such that the new support $\Omega $ satisfies conditions (i)-(iii), namely, (i) any PL map $\beta $ from the ball $B_m(1)$ still keeps $\Omega $ inside an interval $(0,\epsilon _0)$; (ii) the image of $\Omega $ by elements of $B_m(1)$ with different slopes at 0 are disjoint; and (iii) the images of $\Omega $ by elements of $B_m(1)$ with the same slope at 0 are the same. Then we pick up a generating set which involves the elements $\overline{w_1}, \ldots , \overline{w_k}$. By properties (i)-(iii), we again obtain a short relation among $\overline{w_1}, \ldots , \overline{w_k}$ thus a contradiction.  
 
\vspace{1cm}

\section {Girth of subgroups with hyperbolic-like elements} 

\medskip
	
 In this section we will prove Theorem 1.5.
 
 \medskip
 
  Let $d(\Gamma ) = s$, and $m$ be a natural number. We will find $s+2$ generators of $\Gamma $ such that there is no relation of length $m$ or less in $\Gamma $ in these generators. (since $m$ is arbitrary, this proves that $girth(\Gamma ) = \infty $).
 
 \medskip

Let $S = \{X_1, \ldots , X_s\}$ be a finite generating set of $\Gamma $, and $S^{\ast }$ be the symmetrization of $S$, i.e. $S^{\ast } = \{X_1, \ldots , X_s, X_1^{-1}, \ldots , X_s^{-1}\}$. Let also $p_0\in (0,1)$ (one could take $p_0 = \frac{1}{2}$). We can find a natural number $N > 4m$ and a sequence  $0 = c_0 < c_1 < c_2 < \ldots < c_{2N} < c_{2N+1} = 1$ such that the following three conditions are satisfied.

\medskip

 (i) $p_0\in (c_{N}, c_{N+1})$;

\medskip 

 (ii) for all $X\in S^{\star }$, $X(p_0)\notin \{c_1, \ldots , c_{2N}\}$;
 
 \medskip

 (iii) $W(X_1, \ldots , X_s, Y_N)(p_0)\subset (c_1, c_{2N})$ for all reduced words $W$ of the form $W_0Y_N^{n_1}W_1Y_N^{n_2}\ldots W_{k-1}Y_N^{n_k}W_k$ where $Y_N$ is any orientation preserving homeomorphism satisfying the conditions $Y_N(c_i) = c_{i+4}, 1\leq i\leq 2N-4$, $n_1, \ldots , n_k\in \{-1, 1\}$, and $W_i$ is a reduced word in the alphabet $\{X_1^{\pm 1}, \ldots , X_s^{\pm 1}\}$ of length $L_i, 0\leq i\leq k$ where $\displaystyle \sum _{i=0}^kL_i \leq 2m$.

\medskip

 Let $\delta  = \displaystyle \min _{0\leq i\leq 2N} |c_{i+1}-c_i|$ and $\epsilon = \frac{1}{8}\min \{\delta , |p_0-c_{N}|, |p_0-c_{N+1}|\}$. Then we can find a natural number $M > 2m$ and elements $\gamma , \theta \in \Gamma $ such that 

\medskip
 
 (iv) $Fix(\gamma ) = \{0, a_1, a_2, \ldots , a_N, 1\}, Fix(\theta ) = \{0, b_1, b_2, \ldots , b_N, 1\}$ and for all $1\leq i\leq N$, the inequalities $|a_i-c_{2i-1}| < \epsilon $ and $|b_i-c_{2i}| < \epsilon $ hold;

\medskip

 (v) for all $n\geq M$, we have $\gamma ^{\pm n}(U_{\gamma })\subset V_{\gamma }$ where $$U_{\gamma } = \displaystyle \mathop{\sqcup} _{0\leq i\leq N}(a_i + \epsilon , a_{i+1} - \epsilon ), V_{\gamma } = \displaystyle \mathop {\sqcup} _{0\leq i\leq N+1}(a_i - \epsilon , a_i + \epsilon );$$

\medskip

(vi) for all $n\geq M$, we have $\theta ^{\pm n}(U_{\theta })\subset V_{\theta }$ where $$U_{\theta } = \displaystyle \mathop{\sqcup} _{0\leq i\leq N}(b_i + \epsilon , b_{i+1} - \epsilon ), V_{\theta } = \displaystyle \mathop {\sqcup} _{0\leq i\leq N+1}(b_i - \epsilon , b_i + \epsilon ).$$

\bigskip

It is straightforward to make all these arrangements. (Let us also clarify that we define $a_0 = b_0 = 0$ and $a_{N+1} = b_{N+1} = 1$). 

\medskip

Notice that $p_0\in U_{\gamma } \cup U_{\theta }$ and $p_0\notin V_{\gamma } \cup V_{\theta }$. Now, let $r \geq 2M$, and $S' = \{\gamma , \theta , \gamma ^rX_1\theta ^r, \ldots , \gamma ^{sr}X_s\theta ^{sr} \}$.

\medskip

Then if $W_0$ is a non-trivial reduced word in these generators of length at most $m$, and if $W'$ is any suffix of $W_0$ in the alphabet $S'$, then because of (i)-(vi), we have $W'(p_0)\in (c_1, c_N)$. Then clearly (ping-pong argument) $W_0(p_0)\in V_{\gamma }\cup V_{\theta }$ , hence $W_0(p_0)\neq p_0$, hence $W_0\neq 1$. $\square $

\bigskip

 \begin{rem} The assumptions of Theorem 1.5 can be weakened significantly (at the expense of making the statement more technical).
 \end{rem}
 
\medskip

\begin{rem} From the proof we see that one can state a much more general theorem for the girth of groups acting on metric spaces by homeomorphisms. For every non-elementary word hyperbolic group we do have such a theorem indeed (see Theorem 1, [Ak2], where the metric space is the boundary of the group, and for every hyperbolic element we have one attracting and one repelling point). In our case, the metric space is typically non-compact [in the case of Theorem 1.5, the metric space is the non-compact space $(0,1)\cong \mathbb{R}$], and the ``hyperbolic-like" elements have several points (instead of two) which are ``attractive-repelling like" within ``certain compact subspace".  
\end{rem}

 \medskip

 We would like to give a precise definition of a hyperbolic-like element

 \medskip

\begin{defn} Let $X$ be a Hausdorff topological space, $\Gamma $ be a subgroup of {\em Homeo}$(X)$ generated by a finite subset $S\subseteq \Gamma $, $S^{\star } = S\cup S^{-1}\cup \{1\}$, $z\in X, m\in \mathbb{N}$, and $\gamma \in \Gamma $. We say $\gamma $ is $(S,z,m)$-hyperbolic-like if there exists a chain $\Omega _0 \subset \Omega _1\subset \ldots \subset \Omega _m$ of finite subsets of $X$ such that 

 \medskip
 
 (i) $\Omega _0 = \{z\}$;
 
 \medskip
 
 (ii) $s(\Omega _m) \cap \Omega _m = \emptyset , \forall s\in S^{\star }\backslash \{1\}$;
 
 \medskip
 
 (iii) for all $x\in (S^{\star }\backslash \{1\})\Omega _k, 0\leq k\leq m-1$, there exist distinct $p_a, p_r\in \Omega _{k+1}$ such that for all disjoint open neighborhoods $U_{p_a}, U_{p_r}$ of $p_a$ and $p_r$ respectively, there exist $M\in\mathbb{N}$ such that for all $n\geq M$, $\gamma ^{n}(x)\in U_{p_a},  \gamma ^{-n}(x)\in U_{p_r}$.
\end{defn}

  \medskip
	
	 The following theorem is easy to obtain from the proof of Theorem 1.5, we leave the details to the reader.
	
	\medskip
	
	\begin{thm} Let $X$ be a Hausdorff space, $z\in X$, $\Gamma $ be a finitely generated subgroup of {\em Homeo}$(X)$, $S$ be a finite generating set of $\Gamma $. Assume that for all natural $m$ there exists an $(S,z,m)$-hyperbolic-like element of $\Gamma $. Then $girth(\Gamma ) = \infty $. \ $\square $
	\end{thm}
	
			\medskip
			
	\begin{rem} Theorem 4.4. generalizes Theorem 2.1. from [Ak2] which states that any finitely generated non-elementary subgroup of a word hyperbolic group has infinite girth.  
	\end{rem}
	
	\medskip
	
	\begin{rem} It is interesting that the group $F$ (the standard representation of it in $PL_o(I)$) is very rich with hyperbolic-like elements, for the standard finite generating set $S$ of $F$. It is not known to us if the same can be true for a faithful representation of an elementary amenable subgroup of {\em Homeo}$_{+}(\mathbb{R})$.
\end{rem}

 \bigskip

	We will borrow the following definition from [Ak3]
	
	\medskip
	
	\begin{defn} Let $\Gamma $ be a finitely generated group, $d(\Gamma )$ be the minimal number of a generating set of $\Gamma $ and $k\geq d(\Gamma )$ be a positive integer. Then we define $girth _k(\Gamma ) = \displaystyle \sup _{\langle S \rangle = \Gamma , |S| \leq k}girth (\Gamma , S)$. 
	\end{defn}
	
	\medskip
	
	While proving Theorem 1.2, we indeed proved a bit more, namely, for any non-solvable finitely generated group $\Gamma $ of $PL_o(I)$, we proved that $girth _{2d}(\Gamma ) = \infty $ where $d = d(\Gamma )$. Also, in the proof of Theorem 1.5, we indeed prove that $girth _{d+2}(\Gamma ) = \infty $. With a slightly different argument, one can improve this result showing that $girth _{d+1}(\Gamma ) = \infty $; and in Theorem 4.4, one can prove that $girth _{|S|+1}(\Gamma ) = \infty $
	
	\bigskip
	
	{\em Acknowledgment:} I am thankful to Collin Bleak and Matthew Brin for useful discussions. 
	
	\vspace{2cm}
	
	{\bf R e f e r e n c e s}
	
	\vspace{1cm}
	
	[Ak1] Akhmedov, A. \ {\em On the girth of finitely generated groups.} \ Journal of Algebra, {\bf 268}, (2003), 198-208.
	
	\medskip
	
	[Ak2] Akhmedov, A. \ {\em The girth of groups satisfying Tits Alternative.} \ Journal of Algebra, {\bf 287}, (2005), 275-282.
	
	\medskip
	
	[Ak3] Akhmedov, A. \ {\em Quasi-Isometric Rigidity in Group Varieties.} \ Ph.D Thesis. Yale University, 2004.
	
	\medskip
	
	[AST] Akhmedov, A. \ Stein, M. \ Taback, J. \  {\em Free limits of R.Thompson's group $F$.}  \ Geometriae Dedicata, {\bf 155}, (2011), no.1, 163-176.
	
	\medskip
	
	 [Bl] Bleak, C. \ {\em A geometric classification of some solvable groups of homeomorphisms.}, \ Journal of London Mathematical Society (2), {\bf 78}, (2008) no. 2, 352-372.  
	
	 \medskip
	
	[Br1] Brin, M. \ {\em The Free Group of Rank 2 is a Limit of Thompson's Group F.}, \ Groups Geometry Dynamics, {\bf 4} (2010) no.3, 433-454.
	
	\medskip
	
	[Br2] Brin, M. \ {\em Elementary amenable subgroups of R. Thompson's group F.} \ International Journal of Algebra and Computation {\bf 15} (2005) no 4, 619-642.

 \medskip

  [DS] Dranishnikov, A. \ Sapir. M \ {\em On the dimension growth of groups.} \ ArXiv:1008.3868.

 \medskip

  [N1] Nakamura, K.  \ Ph.D Thesis, University of California Davis, 2008.
	
	\medskip
	
	[N2] Nakamura, K. \ The Girth Alternative for Mapping Class Groups. \ to appear in {\em Groups, Geometry, Dynamics.}, http://arxiv.org, arXiv:1105.5422.
	
	\medskip

  [S] Schleimer, S. \ {\em The girth of groups.} \ Preprint.
	
	\medskip
	
	[Y] Yamagata, S. \ {\em The girth of convergence groups and mapping class groups.} \ Osaka Journal of Mathematics, {\bf volume 48}, no.1, 2011.

\end{document}